\documentclass[12pt,reqno]{amsart}
\usepackage{amssymb,amscd,amsbsy}
\usepackage{amsthm}
\newcommand{\lb}{\linebreak}

\newcommand{\vk}{\varkappa}
\newcommand{\z}{\zeta}

\newcommand{\vt}{\vartheta}
\renewcommand{\l}{\lambda}

\newcommand{\f}{\varphi}

\newcommand{\w}{\omega}

\newcommand{\D}{\Delta}
\renewcommand{\L}{\Lambda}

\renewcommand{\O}{\Omega}

\newcommand{\A}{{\mathcal A}}

\newcommand{\1}{{\bf 1}}

\newcommand{\C}{{\Bbb C}}
\newcommand{\T}{{\Bbb T}}
\newcommand{\pp}{{\Bbb P}}
\newcommand{\dd}{{\Bbb D}}
\newcommand{\R}{{\Bbb R}}

\newcommand{\mm}{{\Bbb M}}
\newcommand{\0}{{\boldsymbol{0}}}

\newcommand{\bs}{\boldsymbol}

\newcommand{\rf}[1]{(\ref{#1})}

\newcommand{\df}{\stackrel{\mathrm{def}}{=}}
\newcommand{\dist}{\operatorname{dist}}
\newcommand{\Ker}{\operatorname{Ker}}

\newcommand{\clos}{\operatorname{clos}}

\newcommand{\rank}{\operatorname{rank}}

\newcommand{\eeq}{\end{equation}}
\newcommand{\beq}{\begin{equation}}
\newcommand{\bay}{\begin{eqnarray}}
\newcommand{\ba}{\begin{align*}}
\newcommand{\ea}{\end{align*}}
\newcommand{\ey}{\end{eqnarray}}
\newcommand{\bey}{\begin{eqnarray*}}
\newcommand{\eey}{\end{eqnarray*}}

\newcommand{\be}{\infty}

\newcommand{\bl}{\blacksquare}
\newcommand{\ess}{\operatorname{ess}}

\newcommand{\Range}{\operatorname{Range}}

\newcommand{\Pf}{{\bf Proof. }}

\newcommand{\ov}{\overline}

\newtheorem{thm}{\hspace{\parindent}Theorem}[section]

\newtheorem{cor}[thm]{\hspace{\parindent}Corollary}
\newtheorem{lem}[thm]{\hspace{\parindent}Lemma}

\theoremstyle{remark}

\newtheorem*{rem*}{Remark}

\begin{document}

\newcommand{\vse}{\vspace{.2in}}
\numberwithin{equation}{section}

\title{\bf Analytic approximation of rational matrix functions}
\author{V.V. Peller and V.I. Vasyunin}
\thanks{The first author is partially supported by NSF grant DMS 0200712,
the second author is partially supported by RFBR grant \#05-01-00925}

\maketitle

\newcommand{\mt}{{\mathcal T}}

\begin{abstract}
For a rational matrix function $\Phi$ with poles outside the unit circle,
we estimate the degree of the unique superoptimal approximation
$\A\Phi$ by matrix functions analytic in the unit disk.
We obtain sharp estimates in the case of $2\times2$ matrix functions.
It turns out that ``generically'' $\deg\A\Phi\le\deg\Phi-2$.
We prove that for an arbitrary $2\times2$ rational function $\Phi$,
$\deg\A\Phi\le2\deg\Phi-3$ whenever $\deg\Phi\ge2$.
On the other hand, for $k\ge2$, we construct a $2\times2$ matrix function $\Phi$, for which
$\deg\Phi=k$, while $\deg\A\Phi=2k-3$. Moreover, we conduct a detailed analysis of the situation 
when the inequality $\deg\A\Phi\le\deg\Phi-2$ can violate and obtain best possible results.
\end{abstract}

\section{\bf Introduction}
\setcounter{equation}{0}

\

The problem of uniform
approximation of a given function on the unit circle $\T$ by functions analytic in the
unit disk is a classical problem of approximation theory that is very important in
various applications (see \cite{Pe}). It is well known (follows from a compactness argument)
that for a function $\f$ in $L^\be(\T)$ there is always a best approximant $f\in H^\be$:
$$
\|\f - f\|_\infty = \dist_{L^\infty}(\f, H^\infty).
$$
Though in general a best approximant does not have to be unique, under certain natural
assumptions on $\f$ we do have uniqueness. In particular, uniqueness holds for continuous
functions $\f$; this was established for the first time in \cite{Kh} and then rediscovered
by several other mathematicians.

A function $\f\in L^\be$ on $\T$ is called {\it badly approximable} if the zero function is a
best approximant of $\f$. It is well known (see e.g., \cite{Pe}, Ch.~7, \S\,5) that a nonzero
continuous function is badly approximable if and only if it has constant modulus and negative
winding number with respect to the origin. This criterion can be extended to functions in
$QC\df H^\be+C\bigcap\ov{H^\be+C}$, see \cite{Pe}, Ch.~7, \S\,5.
If $\f$ is a nonzero rational function, then $\f$ is badly approximable if and only if
it has constant modulus and
\bay
\label{+-}
\deg\pp_+\f<\deg\pp_-\f
\ey
(this follows from the above criterion and the argument principle).

It is even more important in applications to consider the problem of analytic approximation of
matrix functions on the unit circle; in systems theory scalar functions correspond
to single input -- single output systems, while matrix-valued functions correspond to multiple
input -- multiple output systems (see e.g., \cite{Fu}, \cite{Fr}, \cite{Pe}).

It turns out, however, that the situation with best analytic approximation of matrix functions
is quite different from the scalar case. Suppose that $\Phi\in L^\be(\mm_{m,n})$, i.e.,
$\Phi$ is a bounded measurable function with values in the space $\mm_{m,n}$ of $m\times n$
matrices. A matrix function $F$ in $H^\be(\mm_{m,n})$ is said to be a {\it best approximant of}
$\Phi$ if
$$
\|\Phi-F\|_{L^\be}=\inf\{\|\Phi-Q\|_{L^\be}:~Q\in H^\be(\mm_{m,n})\}.
$$
Here for a function $\Psi$ in $L^\be(\mm_{m,n})$ we use the notation
$$
\|\Psi\|_{L^\be}\df\ess\sup_{\z\in\T}\|\Psi(\z)\|_{\mm_{m,n}}.
$$

However, in the matrix case even for polynomial matrix functions a best
approximant is almost never unique. Consider the matrix function
$\Phi=\left(\begin{array}{cc}\bar z&\bs{0}\\\bs{0}&\0\end{array}\right)$.
It is easy to see that $\|\Phi\|_{L^\be}=1$. On the other hand,
if $f$ is a scalar function in the unit ball of $H^\be$ then the function
$\left(\begin{array}{cc}\0&\bs{0}\\\bs{0}&f\end{array}\right)\in H^\be\big(\mm_{2,2}\big)$ is
a best approximant of $\Phi$.

In the matrix case it is more natural to impose additional requirements on a best approximant to
be able to choose among best approximants the``very best'' one.
Recall that for a matrix (or an operator on Hilbert space) $A$ the {\it singular value} $s_j(A)$,
$j\ge0$, is, by  definition, the distance from $A$  to the set of matrices (operators) of rank
at most $j$. Clearly, $s_0(A)=\|A\|$.

\medskip

{\bf Definition.}
Given a matrix function $\Phi\in L^\be(\mm_{m,n})$, we define inductively
the sets $\bs{\O}_j$, $0\le j\le\min\{m,n\}-1$, by
$$
\bs{\O}_0=\{F\in H^\be(\mm_{m,n})
:~F~\mbox{minimizes}~\ t_0\df\ess\sup_{\z\in\T}\|\Phi(\z)-F(\z)\|\};
$$
$$
\bs{\O}_j=\{F\in \O_{j-1}:~F~\mbox{minimizes}~\
t_j\df\ess\sup_{\z\in\T}s_j(\Phi(\z)-F(\z))\},\quad j>0.
$$
Functions in $\bigcap\limits_{k\ge 0} \bs{\O}_k = \bs{\O}_{\min\{m,n\}-1} $ are called {\it
superoptimal approximants} of $\Phi$ by bounded analytic matrix functions. The numbers
$t_j=t_j(\Phi)$ are called the {\it superoptimal singular values} of $\Phi$. Note that the set
$\bs{\O}_0$ consists of the best approximants of $\Phi$.

We say that $\Phi$ is {\it badly approximable} if
$$
\|\Phi\|_{L^\be}=\dist_{L^\be}\big(\Phi,H^\be(\mm_{m,n})\big)
$$
and we say that $\Phi$ is
{\it very badly approximable} if the zero function is a superoptimal approximant  of $\Phi$.

\medskip

It was proved in \cite{PY1} (see also \cite{Pe}, Ch.~14, \S\,3) that for a continuous $m\times
n$ matrix function $\Phi$, there exists a unique superoptimal approximant $F$ and
$$
s_j\big(\Phi(\z)-F(\z)\big)=t_j(\Phi),\quad0\le j\le\min\{m,n\}-1,\quad
\mbox{for almost all}\quad\z\in\T.
$$
For a continuous matrix function $\Phi$, {\it we denote by $\A\Phi$ the unique
superoptimal approximant of $\Phi$}.

Note that rational matrix functions play a special role in applications in control theory;
they correspond to transfer functions whose minimal realizations have finite-dimensional
state spaces. Moreover, the dimension of the state space of minimal realization is equal to
the (McMillan) degree $\deg\Phi$ of a matrix rational function $\Phi$, see e.g., \cite{Fu},
\cite{Pe} (see the definition of $\deg\Phi$ in \S\,2). That is why it is important to obtain
sharp estimates of the degree of a rational matrix function.

It is well known (see \cite{Pe}, Ch.~7,
\S\,1) that for a scalar rational function $\f$ the following inequality holds:
$$
\deg\A\f\le\deg\f\quad\mbox{and}\quad\deg\A\f\le\deg\f-1\quad\mbox{unless}\quad\f\in H^\be
$$
(if $\f\in H^\be$, then clearly, $\A\f=\f$).

It was shown in \cite{PY1} (see also \cite{Pe}, Ch. 14, \S\,12) that if $\Phi$ is a rational matrix
function with poles outside $\T$, then $\A\Phi$ is also rational. It is a very important
problem to obtain sharp estimate of the degree of $\A\Phi$ in terms of the degree of $\Phi$.
In particular, it has been an open problem for over 10 years whether \bay \label{za}
\deg\A\Phi\le\deg\Phi. \ey

In this paper we obtain definitive results for $2\times2$ matrix functions. We show in \S\,3
that  \rf{za} holds ``generically''. Moreover, generically in the nondegenerate case
\lb$\rank(\Phi-\A\Phi)=2$ the following stronger inequality holds:
$$
\deg\A\Phi\le\deg\Phi-2.
$$

However, in \cite{PY2} an example of a $2\times2$ rational matrix function $\Phi$ is given,
for which $\rank(\Phi-\A\Phi)=2$, but $\deg\A\Phi\le\deg\Phi-1$.

In this paper we show that actually, the degree of $\A\Phi$ can be greater than the degree of
$\Phi$.  In \S\,4 for $k\ge2$, we give an example of a $2\times2$ matrix function $\Phi$ such
that $\deg\Phi=k$, while $\deg\A\Phi=2k-3$.

This is the worst possible case. We show in \S\,3 that if $k\ge2$ and $\deg\Phi=k$, then
$\deg\A\Phi\le 2k-3$. Note that the last inequality was obtained by state space methods in
\cite{LHG}, but under an additional (rather complicated) assumption on $\Phi$. Note also that the
authors of \cite{LHG} also obtained (under an additional assumption) estimates in the case of
rational matrix functions of an arbitrary size.

Let us describe the results of \S\,3 in more detail. Let
$$
(\pp_+\Phi)(z)\df\sum_{j\ge0}\hat\Phi(j)z^j\quad\mbox{and}\quad\pp_-\Phi\df\Phi-\pp_+\Phi,
$$
where $\hat\Phi(j)$ is the $j$th Fourier coefficient of the restriction of $\Phi$ to $\T$.
Consider the matrix function $\Psi$ defined by
$\Psi=\Phi-\A\Phi$. Then $\Psi$ is a very badly approximable matrix function. Clearly,
$$
\A\pp_-\Psi=\A\pp_-\Phi=\A\Phi-\pp_+\Phi=-\pp_+\Psi
$$
and it is easy to see that
$$
\deg\A\Phi-\deg\Phi\le\deg\pp_+\Psi-\deg\pp_-\Psi
$$
(the degree of a rational matrix function is defined in \S\,2). Thus the problem of estimating
the degree of $\A\Phi$ in terms of the degree of $\Phi$ reduces to the problem of estimating
the degree of $\pp_+\Psi$ in terms of the degree of $\pp_-\Psi$ for very badly approximable
rational matrix functions $\Psi$.

It follows from the results of \cite{PY1} (see also \cite{Pe}, Ch.~14, \S\,5 and \S\,12) that $\Psi$ is a
$2\times2$ very badly approximable rational matrix function if and only if the function
$\Psi$ admits on the unit circle $\T$ a special (thematic) factorization of the form 
\bay
\label{tf} 
\Psi=\left(\begin{array}{cc}\bar w_1&-w_2\\\bar w_2&w_1\end{array}\right)
\left(\begin{array}{cc}t_0u_0&\0\\\0&t_1u_1\end{array}\right) \left(\begin{array}{cc}\bar
v_1&\bar v_2\\-v_2&v_1\end{array}\right), 
\ey 
where $v_1$, $v_2$, $w_1$, and $w_2$ are scalar
rational functions such that $|v_1|^2+|v_2|^2=|w_1|^2+|w_2|^2$ on $\T$, $v_1$ and $v_2$ have
no common zeros in $\dd$, and $w_1$ and $w_2$ have no common zeros in $\dd$; $t_0\ge
t_1\ge0$; and $u_0$ and $u_1$ are scalar badly approximable rational functions such that
$|u_0|=|u_1|=1$ on $\T$.

We prove in \S\,3 that if $\Psi$ is nondegenerate (i.e., $t_1>0$) and the inequality 
\bay
\label{ob} 
\deg\pp_+\Psi\le\deg\pp_-\Psi-2 
\ey 
is violated, then there is a point $\l$ in the
unit disk $\dd$ such that $u_1$ has a pole at $\l$, $u_0$ has a zero at $\l$, and both
$\bs{v}=\left(\begin{array}{c}v_1\\v_2\end{array}\right)$ and
$\bs{w}=\left(\begin{array}{c}w_1\\w_2\end{array}\right)$ have poles at $1/\bar\l$. We also
show in \S\,3 that for fixed $t_0, u_0, u_1, v_1, v_2, w_1$, and $w_2$, inequality \rf{ob} can
be violated for no more than $\deg\pp_-u_1$ numbers $t_1\in(0,t_0]$. Moreover, we show
in \S\,3 that for $\Psi$ given by \rf{tf}, the sum of $\left[\deg\pp_+\Psi+2-\deg\pp_-\Psi\right]_+$
over all numbers $t_1\in(0,t_0]$ cannot be greater than $\deg\pp_-u_1$. Here for $x\in\R$, we use the notation
$$
[x]_+=\max\{x,0\}.
$$

In \S\,4 we construct an example of a matrix function of the form \rf{tf} such that inequality
\rf{ob} is violated for precisely $\deg\pp_-u_1$ values of $t_1\in(0,t_0]$. This example 
also shows that that the above estimate obtained in \S\,3 is best possible.

In \S\,2 we collect information on Hankel and Toeplitz operators, degrees of rational matrix
functions, and very badly approximable matrix functions.

\

\section{\bf Preliminaries}
\setcounter{equation}{0}

\

{\bf 2.1. Hankel and Toeplitz operators on spaces of vector functions.} For a matrix
function $\Phi$ in $L^\be(\mm_{m,n})$, the {\it Hankel operator}
$$
H_\Phi:H^2(\C^n)\to H^2_-(\C^m)\df L^2(\C^m)\ominus H^2(\C^m),
$$
is defined by
$$
H_\Phi f=\pp_-(\Phi f),
$$
where $\pp_-$ is the orthogonal projection from $L^2(\C^m)$ onto $H^2_-(\C^m)$. Here
$H^2(\C^n)$ is the Hardy space of $\C^n$-valued functions that is identified naturally with a
subspace of $L^2(\C^n)$.

It is well known that
$$
\|H_\Phi\|=\dist_{L^\be}\big(\Phi,H^\be(\mm_{m,n})\big).
$$
The operator $H_\Phi$ is compact if and only if $\Phi\in\big(H^\be+C\big)(\mm_{m,n})$.

For $\Phi\in L^\be(\mm_{m,n})$, the {\it Toeplitz operator}
$$
T_\Phi:H^2(\C^n)\to H^2(\C^m),
$$
is defined by
$$
T_\Phi f=\pp_+(\Phi f),
$$
where $\pp_+$ is the orthogonal projection from $L^2(\C^m)$ onto $H^2(\C^m)$.

If $\Phi$ is a continuous square matrix function, then $T_\Phi$ is {\it Fredholm}
(an operator $T$ on Hilbert space is called Fredholm if its range is closed and
both $\Ker T$ and $\Ker T^*$ are finite-dimensional) if and only if $\det\Phi$ does not vanish
on $\T$.

We refer the reader to \cite{Pe} for more detailed information on Hankel and Toeplitz operators
on spaces of vector functions.

\

{\bf 2.2. The degree a of rational matrix function.} To define the (McMillan) degree of a
rational  matrix function, we define first Blaschke--Potapov products. An $n\times n$ matrix
function
$B$ is called a  {\it finite Blaschke--Potapov product} if it admits a factorization
\beq
\label{BP}
B=UB_1B_2\cdots B_k,
\end{equation}
where
\[
B_j(z)=\frac{\l_j-z}{1-\bar{\l}_jz}P_j+(I-P_j)
\]
for some $\l_j$ in $\dd$, orthogonal projections $P_j$ on $\C^n$,
and a unitary matrix $U$.
The degree of the Blaschke--Potapov product (\ref{BP}) is defined by
\[
\deg B=\sum^k_{j=1}\rank P_j.
\]

Let now $\Phi$ be an $m\times n$ rational function without poles in
$\hat\C\setminus\clos\dd$. We say that $\Phi$ has {\it \textup(McMillan\textup) degree} $d>0$
(in other words, $\deg\Phi=d$) if there exists an $n\times n$ Blaschke--Potapov product $B$
of degree $d$ such that $\Phi B\in H^\be(\mm_{m,n})$, but there is no such Blaschke--Potapov
product of degree $d-1$. If $\Phi$ is a constant function, we put $\deg\Phi=0$. It is well
known (see e.g., \cite{Pe}, Ch.~2, \S\,5) that $\deg\Phi$ is the rank of the Hankel operator
$H_\Phi:H^2(\C^n)\to H^2_-(\C^m)$ for rational functions $\Phi$ with poles in
$\dd$ and 
\bay
\label{tr}
\deg\Phi=\deg\Phi^{\rm t}
\ey 
for an arbitrary rational matrix function $\Phi$.

Note that it follows from \rf{tr} that $\deg\Phi=d>0$ if and only if
there exists a Blaschke--Potapov product $B$
of degree $d$ such that $B\Phi \in H^\be(\mm_{m,n})$, but there is no such Blaschke--Potapov
product of degree $d-1$.

If $\Phi$ is a rational matrix function without poles in $\clos\dd$, we say that
$\deg\Phi=\deg\Phi_\flat$, where the rational matrix function $\Phi_\flat$ is defined by
$\Phi_\flat(z)=\Phi(1/z)$.

For a matrix function $\Phi$ without poles on the unit circle $\T$ we put
$$
\deg\Phi=\deg\pp_-\Phi+\deg\pp_+\Phi.
$$

It is easy to see that if for $r>0$, the function $\Phi_r$ is defined by
$$
\Phi_r(z)=\Phi(rz)
$$
and both $\Phi$ and $\Phi_r$ have no poles on $\T$, then $\deg\Phi=\deg\Phi_r$.
If $\Phi$ is an arbitrary rational matrix function, we define the {\it degree} of $\Phi$ as
$$
\deg\Phi=\deg\Phi_r,
$$
where $r$ is a positive number such that $\Phi_r$ has no poles on $\T$.

For a rational matrix function $\Phi$ and $\l\in\hat{\C}$ we define the {\it
\textup(McMillan\textup) degree of \,$\Phi$ at $\l$} as the degree of the principal part of
$\Psi$ at $\l$ and we use the notation $\deg_\l\Phi$ for the degree of $\Phi$ at $\l$.

Clearly,
$$
\deg\Phi=\sum_{\l\in\hat\C}\deg_\l\Phi.
$$

For a subset $\L$ of $\hat{\C}$, we use the notation
$$
\deg_\L\Phi=\sum_{\l\in\L}\deg_\l\Phi.
$$

It is easy to see that for $\L\subset\dd$, $\deg_\L\Phi=\vk>0$ if and only if there exists 
a Blaschke--Potapov product $B$
of degree $\vk$ such that $B\Phi$ has no pole in $\L$, but there is no such Blaschke--Potapov
product of degree $\vk-1$.

We need the following elementary facts about the degrees of rational functions.

\begin{lem}
\label{dt}
Let $\Phi$ be a square rational matrix function. Then
$$
\deg_\L\det\Phi\le\deg_\L\Phi
$$
for every $\L\subset\hat\C$.
\end{lem}

\Pf Clearly, it suffices to consider the case when 
$\L\subset\dd$. Let $d=\deg_\L\Phi$ and let $B$ be a Blaschke--Potapov product of degree $d$ such that
$\Phi B$ has no poles in $\L$. Then
$\det\Phi B=\det\Phi\det B$ has no poles in $\L$.
It follows immediately from the definition of Blaschke--Potapov products that $\det B$ is a
scalar Blaschke product of degree $d$. Thus $\deg_\L\det\Phi\le d$. \hfill$\bl$

\medskip

Suppose now that
$\Phi=\left(\begin{array}{cc}\f_{11}&\f_{12}\\\f_{21}&\f_{22}\end{array}\right)$  is a
$2\times2$ matrix function. We denote by $\Phi^{\rm ad}$ the adjugate matrix function:
$$
\Phi^{\rm ad}\df\left(\begin{array}{rr}\f_{22}&-\f_{21}\\-\f_{12}&\f_{11}\end{array}\right).
$$

\begin{lem}
\label{jug}
Let $\Phi$ be a $2\times2$ rational matrix function. Then
$$
\deg\Phi^{\rm ad}=\deg\Phi.
$$
\end{lem}

\Pf For the unitary matrix $J=\begin{pmatrix} 0&-1\\1&\phantom{-}0\end{pmatrix}$ we
have $\Phi^{\rm ad}=J\Phi J^*$. It is easy to see from the definition of
the degree that the degree does not change when we multiply a matrix function by
a constant unitary matrix. \hfill$\bl$

\begin{lem}
\label{polu}
Suppose that $\Phi_1$ is an $m\times n$ rational matrix function and $\Phi_2$ is an 
$n\times k$ rational matrix function. Then
$$
\deg_\L(\Phi_1\Phi_2)\le\deg_\L\Phi_1+\deg_\L\Phi_2
$$
for every $\L\subset\hat\C$.
\end{lem}

\Pf Clearly, it suffices to prove this inequality for $\L\subset\dd$.
Let us show that if $\Psi$ is a rational function without poles on $\T$ and $B$ is a finite Blaschke--Potapov product,
then
\bay
\label{-B}
\deg_\L(\Psi B)\ge\deg_\L\Psi-\deg B.
\ey
Let $B_\star$ be a Blaschke--Potapov product of degree $\deg_\L(\Psi B)$
such that $\Psi BB_\star$ has no pole in $\L$. Then $\deg(BB_\star)\ge\deg_\L\Psi$. Together with the obvious
equality $\deg(BB_\star)=\deg B+\deg B_\star$ this proves \rf{-B}.

Similarly, if $B$ is a finite Blaschke--Potapov product, then
\bay
\label{-Bl}
\deg_\L(B\Psi)\ge\deg_\L\Psi-\deg B.
\ey

Let now $B_1$ be a Blaschke--Potapov product of degree $\deg_\L\Phi_1$ such that 
$B_1\Phi_1$ has no pole in $\L$ and let $B_2$ be a Blaschke--Potapov product of degree $\deg_\L\Phi_2$ such that 
$\Phi B_2$ has no pole at $\L$. Then by \rf{-B} and \rf{-Bl},
$$
0=\deg_\L(B_1\Phi_1\Phi_2B_2)\ge\deg_\L(B_1\Phi_1\Phi_2)-\deg B_2\ge\deg_\L(\Phi_1\Phi_2)-\deg B_1-\deg B_2
$$
which proves the result. \hfill$\bl$

\medskip

{\bf Remark.} In fact, Lemma \ref{dt} in the special case of $2\times2$ functions is an easy consequence
of Lemma \ref{jug}:
$$
(\deg_\l\det\Phi)^2=\deg_\l((\det\Phi)I_{2\times2})=\deg_\l(\Phi^{\rm ad}\Phi)\le
\deg_\l\Phi^{\rm ad}\deg_\l\Phi=(\deg_\l\Phi)^2.
$$

\

{\bf2.3. Thematic factorizations and very badly approximable matrix functions.} For
simplicity, we state here results on badly approximable matrix functions in the $2\times2$
case. Note that all results mentioned in this subsection are also valid in the case of matrix
functions of an arbitrary size;  we refer the reader to \cite{Pe} (Ch.~14)
for the general case. However, in the case of an arbitrary size
thematic factorizations look more complicated.

In \cite{PY1} very badly approximable matrix functions were characterized in terms of certain
special factorizations ({\it thematic factorizations}). In particular, it follows from the
results of \cite{PY1} that a continuous matrix function $\Psi$ on $\T$ is very badly
approximable if and only if it admits a factorization 
\bay 
\label{si}
\Psi=\left(\begin{array}{cc}\bar w_1&-w_2\\\bar w_2&w_1\end{array}\right)
\left(\begin{array}{cc}t_0u_0&\0\\\0&t_1u_1\end{array}\right) \left(\begin{array}{cc}\bar
v_1&\bar v_2\\-v_2&v_1\end{array}\right), 
\ey 
where $0\le t_1\le t_0$, $u_0$ and $u_1$ are
badly approximable unimodular (i.e., $|u_0|=|u_1|=1$ almost everywhere on
$\T$) functions of class $QC= H^\be+C\bigcap\ov{H^\be+C}$;
$\bs{v}=\left(\begin{array}{c}v_1\\v_2\end{array}\right)$ and
$\bs{w}=\left(\begin{array}{c}w_1\\w_2\end{array}\right)$ are {\it inner and co-outer column
functions} (i.e., $v_1$, $v_2$, $w_1$, and $w_2$ are $H^\be$ functions satisfying
$|v_1|^2+|v_2|^2=|w_1|^2+|w_2|^2=1$ almost everywhere on $\T$ and such that the functions
$v_1$ and $v_2$ are coprime and the functions $w_1$ and $w_2$ are
coprime\footnote {Recall that $H^\be$ functions are called coprime if they have no common
nonconstant inner divisor.}). Moreover, in this case $t_0=t_0(\Psi)$ and $t_1=t_1(\Psi)$ are
the superoptimal singular values of $\Psi$. We refer the reader to \cite{Pe}, Ch.~14, \S\,5.
Note also that another (more geometric) description of very badly approximable matrix
functions is given in \cite{PT}.

It was also shown in \cite{PY1} that in case $\Psi$ is nondegenerate (i.e., $\rank\Psi=2$
almost everywhere on $\T$), the Toeplitz operator $T_{z\Psi}:H^2(\C^2)\to H^2(\C^2)$ has
dense range (see also \cite{Pe}, Ch.~14, Theorem 5.4).

It was also shown in \cite{PY1} that if $\Psi$ is rational, then all functions in
the factorization \rf{si} are rational. Moreover, the following result holds.

\begin{lem}
\label{kn}
The following inequalities hold:
\bay
\label{nu}
\deg\pp_-u_0\le\deg\pp_-\Psi,
\ey
\bay
\label{nv}
\deg\bs{v}\le\deg\pp_-\Psi-1,
\ey
and
\bay
\label{nw}
\quad\deg\bs{w}\le\deg\pp_-\Psi-1.
\ey
\end{lem}

The proof of Lemma \ref{kn} is implicitly contained in the proof of Lemma 12.2 of Ch.~14 of\cite{Pe}.
We give here an explicit proof.

\medskip

{\bf Proof of Lemma \ref{kn}.} Multiplying~\rf{si} on the left by $\bs{w}^{\rm t}$, we
arrive at the following identity:
\bay
\label{uvwP}
\bs{w}^{\rm t}\Psi=t_0u_0\bs{v}^*.
\ey
Applying this equality to the vector function $\bs{v}$, we obtain
$$
\bs{w}^{\rm t}\Psi\bs{v}=t_0u_0.
$$
Let $\l\in\dd$. Since $\bs{v}$ and $\bs{w}^{\rm t}$ are analytic in $\dd$, we have
$\deg_\l\bs{w}^{\rm t}=\deg_\l\bs{v}=0$, and so by Lemma~\ref{polu},
\begin{align}
\label{sa}
\deg_\l u_0=\deg_\l(\bs{w}^{\rm t}\Psi\bs{v})\le
\deg_\l\bs{w}^{\rm t}+\deg_\l\Psi+\deg_\l\bs{v}=\deg_\l\Psi.
\end{align}
To obtain \rf{nu}, one has to take the sum in \rf{sa} over $\dd$.

To prove \rf{nv}, we multiply \rf{uvwP} by $1/u_0$:
$$
\frac1{u_0}\bs{w}^{\rm t}\Psi=t_0\bs{v}^*,
$$
whence
$$
\deg\bs{v}=\deg\bs{v}^*=\sum_{\l\in\dd}\deg_\l\frac1{u_0}\bs{w}^{\rm t}\Psi=
\sum_{\l\in\dd}\deg_\l\frac1{u_0}\Psi.
$$
If $\l$ is a pole of $\Psi$, then by \rf{sa}, $\deg_\l\Psi\ge\deg_\l u_0$, and since $u_0$ is
a scalar function, we have
$$
\deg_\l\frac1{u_0}\Psi\le\deg_\l\Psi-\deg_\l u_0.
$$
If $\l$ is a regular point of $\Psi$, then
$$
\deg_\l\frac1{u_0}\Psi=\deg_\l\frac1{u_0}=\deg_{1/\ov{\l}}u_0.
$$
Hence, for arbitrary $\l\in\dd$, we have
$$
\deg_\l\frac1{u_0}\Psi\le\deg_\l\Psi-\deg_\l u_0+\deg_{1/\ov{\l}}u_0.
$$
Summation over $\l\in\dd$ yields
$$
\deg\bs{v}\le\deg\pp_-\Psi-\deg\pp_-u_0+\deg\pp_+u_0.
$$
Inequality in \rf{nv} follows now from \rf{+-}. To prove \rf{nw}, one has to apply \rf{nv} to
the matrix function $\Psi^{\rm t}$. \hfill$\bl$

\medskip

To conclude this subsection, we state one more inequality:
$$
\deg\pp_-u_1\le\deg\pp_-\Psi-\mu,
$$
where $\mu$ is the multiplicity of the singular value $1$ of the Hankel operator $H_{u_0}$.
This is a consequence of Theorem 8.1 of Ch.~14 of \cite{Pe}. Thus
\bay
\label{pu1}
\deg\pp_-u_1\le\deg\pp_-\Psi-1.
\ey

\

{\bf2.4. Singular values of $\bs{H_U}$ and $\bs{H_{U^*}}$ for unitary-valued matrix functions
$\bs{U}$.} Let $U$ be a unitary-valued $n\times n$ matrix function such that the Toeplitz
operator $T_U:H^2(\C^n)\to H^2(\C^n)$ has dense range. Then the operator $H^*_{U^*}H_{U^*}$
is unitarily equivalent to the restriction of $H^*_UH_U$ to $\Ker T_U$ (see \cite{Pe}, Ch.~4,
\S\,4). For scalar unimodular functions this result was proved in \cite{PK}. It is easy to see
that
$$
\Ker T_U=\{f\in H^2(\C^n):~\|H_Uf\|=\|f\|\}.
$$

In particular, if $H_U$ is compact, then $\dim\Ker T_U<\be$, and so \bay \label{sj}
s_j(H_{U^*})=s_{j+\mu}(H_U),\quad\mbox{where}\quad \mu=\dim\Ker T_U. \ey In other words,
$\mu$ is the multiplicity of the singular value 1 of $H_U$. This means that to obtain 
the singular values of $H_{U^*}$, one has to remove from the singular
values of $H_U$ all singular values equal to 1.

Suppose now that $\Psi$ is a continuous very badly approximable function of size $n\times n$ such that
$t_j(\Psi)=1$, $0\le j\le n-1$. Then $\Psi$ is unitary-valued and the Toeplitz operator
$T_{z\Psi}$ has dense range (see \S\,2.3). Since $T_{z\Psi}$ is Fredholm (see \S\,2.1), $\Range
T_{z\Psi}=H^2(\C^n)$. In particular, for any constant vector $c$ there
exists a function $f\in H^2(\C^n)$ such that $T_{z\Psi}f=c$. It follows that $f\in\Ker
T_\Psi$, and so $\dim\Ker T_\Psi\ge n$. Thus it follows from \rf{sj} that
$$
\rank H_{\Psi^*}\le \rank H_\Psi-n.
$$
In other words
\bay
\label{-n}
\deg\pp_+\Psi\le\deg\pp_-\Psi-n.
\ey

\

\section{\bf Upper Estimates}
\setcounter{equation}{0}

\

In this section we are going to obtain sharp estimates of the degree of the superoptimal
approximant $\A\Phi$ of a rational $2\times2$ matrix function $\Phi$ in terms of the degree
of $\Phi$. As we have observed in \S\,1, this problem reduces to the problem for a very badly
approximable rational matrix function $\Psi$, to estimate the degree of $\pp_+\Psi$ in terms
of the degree of $\pp_-\Psi$.

Let $\Psi$ be a very badly approximable rational matrix function with $t_1(\Psi)\ne0$.
Without loss of generality we may assume that $t_0(\Psi)=\|H_\Psi\|=1$. Then $\Psi$ admits a
thematic factorization
\bay
\label{tem}
\Psi=\left(\begin{array}{cc}\bs{w}^\#&\xi\end{array}\right)
\left(\begin{array}{cc}u_0&\0\\\0&t_1(\Psi)u_1\end{array}\right)
\left(\begin{array}{c}\bs{v}^*\\\vt^{\rm t}\end{array}\right)
\ey
Here $u_0$ and $u_1$ are
badly approximable scalar rational functions,
$\bs{v}=\left(\begin{array}{c}v_1\\v_2\end{array}\right)$ and
$\bs{w}=\left(\begin{array}{c}w_1\\w_2\end{array}\right)$ are rational inner and co-outer
column functions, and
\bay
\label{dop}
\xi=\left(\begin{array}{c}-w_2\\w_1\end{array}\right)\quad\mbox{and}\quad
\vt=\left(\begin{array}{c}-v_2\\v_1\end{array}\right)
\ey
(see \S\,2.3).

For a rational matrix function $G$ we are going to use the notation
$$
G^\#(\z)\df\ov{G(1/\bar\z)}
$$
and
$$
G^*=\big(G^\#\big)^{\rm t}.
$$
Clearly, $G^\#$ and $G^*$ are rational matrix functions and
$$
G^\#(\z)=\ov{G(\z)}\quad\mbox{for}\quad\z\in\T.
$$

We consider the parametric family $\{\Psi_{[t]}\}_{t>0}$ of very badly approximable matrix
functions:
\bay
\label{Phit}
\Psi_{[t]}=\left(\begin{array}{cc}\bs{w}^\#&\xi\end{array}\right)
\left(\begin{array}{cc}u_0&\0\\\0&tu_1\end{array}\right)
\left(\begin{array}{c}\bs{v}^*\\\vt^{\rm t}\end{array}\right).
\ey
As we have discussed in
\S\,2.3, $\Psi_{[t]}$ is very badly approximable for $t\le1$.

In this section we prove that for all but finitely many points $t\in(0,1]$ we have
\bay
\label{nd}
\deg\pp_+\Psi_{[t]}\le\deg\pp_-\Psi_{[t]}-2
\ey
and for all $t\in(0,1]$ we have
$$
\deg\pp_+\Psi_{[t]}\le2\deg\pp_-\Psi_{[t]}-3.
$$

In the case $t_0(\Psi)\ne0$ and $t_1(\Psi)=0$ we prove that
$$
\deg\pp_+\Psi\le\deg\pp_-\Psi-1.
$$

Let $\Psi=\Psi_{[t]}$ be a rational very badly approximable matrix function such that
\lb$t_1(\Psi)\ne0$. By \rf{tem},
\bay
\label{Phi}
\Psi_{[t]}
=u_0\bs{w}^\#\bs{v}^*+tu_1\xi\vt^{\rm t}.
\ey

Obviously, \bay \label{T} \det\left(\begin{array}{cc}\bs{w}^\#&\xi\end{array}\right)(\z)
=\det\left(\begin{array}{c}\bs{v}^*\\\vt^{\rm t}\end{array}\right)(\z)=1 \ey for all
$\z\in\T$. Since the matrix functions in \rf{T} are rational, it follows that \rf{T} holds
for all $\z\in\C$. Thus we have 
\bay 
\label{det}
\det\Psi_{[t]}(\z)=tu_0(\z)u_1(\z),\quad\z\in\C. 
\ey

We start with the inequality
\bay
\label{t=1}
\deg\pp_+\Psi_{[1]}\le\deg\pp_-\Psi_{[1]}-2
\ey
(see \rf{-n}).

Let $t\in(0,1]$. Suppose that both functions $u_0\bs{w}^\#\bs{v}^*$ and $tu_1\xi\vt^{\rm t}$
in \rf{Phi} have a pole at $\l\in\hat\C$. We say that we have a {\it cancellation}
at $\l$ on the right-hand side of \rf{Phi} if
$$
\deg_\l(u_0\bs{w}^\#\bs{v}^*+tu_1\xi\vt^{\rm t})<
\max\{\deg_\l u_0\bs{w}^\#\bs{v}^*,\deg_\l u_1\xi\vt^{\rm t}\}.
$$
Since both matrix functions $u_0\bs{w}^\#\bs{v}^*$ and $u_1\xi\vt^{\rm t}$ have rank one,
for a given $\l$, we can have a cancellation at $\l$ for no more than one $t$.

\medskip

{\bf Definition.} We say that $t$ is a {\it nondisturbing number for the family}
$\{\Phi_{[t]}\}$ (or simply a {\it nondisturbing number}) if for this $t$ there is no
cancellation of poles on the right-hand side of \rf{Phi}. We say that $t$ is a {\it
nondisturbing number for a set $\L\subset\hat\C$} if there is no cancellation of poles for $\l\in\L$.
Otherwise we say that $t$ is a
{\it disturbing number} ({\it disturbing number for $\L$}\,).

\medskip

Note that the function $\bs{w}^\#\bs{v}^*$ is analytic in $\hat\C\setminus\clos\dd$, and so
if $u_0\bs{w}^\#\bs{v}^*$ has a pole at $\l\in\hat\C\setminus\clos\dd$, then $\l$ is a pole
of $u_0$. Similarly, the function $\xi\vt^{\rm t}$ is analytic in $\dd$, and so if
$tu_1\xi\vt^{\rm t}$ has a pole at $\l\in\dd$, then $u_1$ has a pole at $\l$. Thus there can
be at most $\deg\pp_+u_0$ disturbing numbers $t$ for the exterior of the unit disk and at
most $\deg\pp_-u_1$ disturbing numbers $t$ for the unit disk. Recall that
\bay
\label{ner}
\deg\pp_+u_0\le\deg\pp_-u_0-1\le\deg\pp_-\Psi_{[t]}-1\quad\mbox{for all}\quad t\in(0,1]
\ey
and
\bay
\label{D}
\deg\pp_-u_1\le\deg\pp_-\Psi_{[t]}-1\quad\mbox{for all}\quad t\in(0,1]
\ey
(see \rf{+-}, \rf{nu}, and \rf{pu1}).

We are going to prove that inequality \rf{nd} holds for all nondisturbing numbers $t$ for the
unit disk. If 1 is a nondisturbing number for the exterior of the unit disk, then \rf{nd} is
an immediate consequence of \rf{t=1}. In the case when 1 is a disturbing number for the
exterior of the unit disk, we prove that 1 must also be a disturbing point for the unit disk
and the number of cancelled poles (counted with multiplicities) outside the unit disk is at
least the number of cancelled poles in the unit disk. Clearly, this would imply \rf{nd} for
all $t$ nondisturbing for the unit disk.

We have from \rf{Phi}
\bay
\label{xi}
\xi^*\Psi_{[t]}=tu_1\vt^{\rm t}
\ey
and
\bay
\label{vt}
\Psi_{[t]}\vt^\#=tu_1\xi.
\ey
Multiplying \rf{xi} by the adjugate matrix $\Psi_{[t]}^{\rm ad}$ on the right and multiplying
\rf{vt} by
$\Psi_{[t]}^{\rm ad}$ on the left, we obtain
$$
\big(\det\Psi_{[t]}\big)\xi^*=tu_1\vt^{\rm t}\Psi_{[t]}^{\rm ad}\quad\mbox{and}
\quad\big(\det\Psi_{[t]}\big)\vt^\#=tu_1\Psi_{[t]}^{\rm ad}\xi.
$$
It follows now from \rf{det} that
\bay
\label{ad}
u_0\xi^*=\vt^{\rm t}\Psi_{[t]}^{\rm ad}\quad\mbox{and}
\quad u_0\vt^\#=\Psi_{[t]}^{\rm ad}\xi.
\ey

\begin{lem}
\label{krat} Let $\Psi$ be a matrix function of the form {\em\rf{tem}} such that
$t_1=t_1(\Psi)\ne0$. Suppose that $u_0$ has a pole of multiplicity $d$ at $\l$\textup,
$|\l|>1$\textup, and $\deg_\l\Psi=d-l$\textup, $l>0$. Then $u_1$ has a zero at $\l$ of
multiplicity at least $l$\textup, and both $\bs{v}$ and $\bs{w}$ have poles at $\l$ of
multiplicity at least $l$.
\end{lem}

\Pf By Lemma \ref{dt},
$$
\deg_\l\det\Psi\le\deg_\l\Psi=d-l.
$$
By \rf{det}, $u_1$ has a zero at $\l$ of multiplicity at least $l$. Note that $\xi^*(\l)\ne0$
and $\vt^\#(\l)\ne0$ (this follows from the fact that $\bs{v}$ and $\bs{w}$ are co-outer).
Thus $\deg_\l u_0\xi^*=\deg_\l u_0\vt^\#=\deg_\l u_0=d$. On the other hand, by
Lemma~\ref{jug}, $\deg_\l\Psi^{\rm ad}=\deg_\l\Psi=d-l$.
Hence, by \rf{ad}, $\deg_\l\xi\ge l$ and $\deg_\l\vt^{\rm t}\ge
l$, i.e., both $\bs{v}$ and $\bs{w}$ have poles at
$\l$ of multiplicity at least $l$. \hfill$\bl$

\medskip

Put
$$
J=\left(\begin{array}{cc}0&-1\\1&0\end{array}\right)
$$
and
$$
A=\bs{w}\bs{v}^{\rm t}.
$$
We have
\bay
\label{J}
\Psi_{[t]}=u_0A^\#-tu_1JAJ.
\ey

\

\begin{lem}
\label{pog}
Let $|\l|>1$. Suppose that for some $t<1$,
\bay
\label{nar}
\deg_\l\Psi_{[t]}-\deg_\l\Psi_{[1]}=l>0.
\ey
Then for all $t<1$,
\bay
\label{vn}
\deg_{1/\bar\l}\Psi_{[t]}-\deg_{1/\bar\l}\Psi_{[1]}=l>0.
\ey
\end{lem}

\Pf We have already observed that since both matrix functions $u_0\bs{w}^\#\bs{v}^*$
and $u_1\xi\vt^{\rm t}$ have rank one, it follows from \rf{Phi} that if \rf{nar} holds for some
$t\in(0,1)$, then it holds for all $t\in(0,1)$. The same is true for \rf{vn}.

Let $d$ be the multiplicity of the pole of $u_0$ at $\l$. Since $1$ is a disturbing number
for the point $\l$, we have
$$
\deg_\l(u_0\bs{w}^\#\bs{v}^*)=\deg_\l(u_1\xi\vt^{\rm t}).
$$
Note that the vector functions $\bs{w}^\#$, $\bs{v}^*$ are analytic and nonzero for $|\l|>1$. Therefore
$$
\deg_\l\Psi_{[t]}=\deg_\l(u_0\bs{w}^\#\bs{v}^*)=\deg_\l u_0=d\qquad\text{for all}\qquad t\ne1
$$
and at the disturbing number the degree drops by $l$: $\deg_\l\Psi_{[1]}=d-l$.

Since for a rational function $u$ unimodular on $\T$, we have $u^\#=1/u$, it follows from
\rf{J} that
\bay
\label{1t}
\Psi_{[t]}^\#=u_0^\#A-tu_1^\#JA^\#J=\frac1{u_0u_1}\left(u_1A-tu_0JA^\#J\right)=
-\frac{t}{u_0u_1}J\Psi_{[1/t]}J.
\ey

By Lemma \ref{krat}, the function $\frac1{u_0u_1}$ at $\l$ cannot have a zero of multiplicity
greater than $d-l$. Thus it follows from \rf{1t} that $1$ is a disturbing number for  the family $\{\Psi^\#_{[t]}\}$
for the point $\l$ and the multiplicity of the pole of $\Psi^\#_{[1]}$ at
$\l$ must also drops by $l$. This means that $t=1$ is a disturbing
number for the family $\{\Psi_{[t]}\}$ for the point $1/\bar\l$ and the multiplicity of the pole of $\Psi_{[1]}$ at
$1/\bar\l$ must also drop by $l$.\hfill$\bl$

\begin{thm}
\label{gen}
For all nondisturbing numbers $t\in(0,1]$ for the unit disk inequality
{\em\rf{nd}} holds. There can be at most $\deg\pp_-u_1$ numbers $t\in(0,1)$, for which
inequality {\em\rf{nd}} can be violated.
\end{thm}

\Pf Recall that for all nondisturbing numbers $t$, $\deg_\L\Psi_{[t]}$
remains the same. Let us first consider the
case when $1$ is a nondisturbing number for the family $\{\Psi_{[t]}\}$ for 
the exterior of the unit disk. Then for all nondisturbing numbers $t$
for the unit disk we have
\begin{align*}
\deg_{\hat\C\setminus\dd}\Psi_{[t]}&=\deg\pp_+\Psi_{[t]}\le\deg\pp_+\Psi_{[1]},
\\
\deg_{\dd}\Psi_{[t]}&=\deg\pp_-\Psi_{[t]}\ge\deg\pp_-\Psi_{[1]}.
\end{align*}
Therefore inequality~\rf{t=1} yields
$$
\deg\pp_-\Psi_{[t]}-\deg\pp_+\Psi_{[t]}\ge\deg\pp_-\Psi_{[1]}-\deg\pp_+\Psi_{[1]}\ge2
$$
for all nondisturbing numbers $t$ for the unit disk.

In the case when $1$ is a disturbing number for the family $\{\Psi_{[t]}\}$ for some $\l$, $|\l|>1$, we use
Lemma~\ref{pog}. Let $\L$ be the set of all points $\l$ satisfying the hypotheses of
Lemma~\ref{pog}, and put $\L_\#=\{1/\bar\l\colon\l\in\L\}$. Taking the sum in \rf{nar} and \rf{vn}
over $\L$ yields
$$
\deg_{\L_\#}\Psi_{[t]}-\deg_\L\Psi_{[t]}=\deg_{\L_\#}\Psi_{[1]}-\deg_\L\Psi_{[1]}.
$$
For all $\l\in\hat\C\setminus(\dd\cup\L)$ the number $1$ is nondisturbing, whence
$$
\deg_{\hat\C\setminus(\dd\cup\L)}\Psi_{[t]}\le\deg_{\hat\C\setminus(\dd\cup\L)}\Psi_{[1]}.
$$
And finally, if $t$ is a nondisturbing number for the unit disc, then we have
\begin{align*}
\deg\pp_-\Psi_{[t]}-\deg\pp_+\Psi_{[t]}\!&=
(\deg_{\dd\setminus\L_\#}\Psi_{[t]}-\deg_{\hat\C\setminus(\dd\cup\L)}\Psi_{[t]})+
(\deg_{\L_\#}\!\Psi_{[t]}-\deg_\L\Psi_{[t]})
\\
&\ge (\deg_{\dd\setminus\L_\#}\!\Psi_{[1]}-\deg_{\hat\C\setminus(\dd\cup\L)}\!\Psi_{[1]})+
(\deg_{\L_\#}\!\Psi_{[1]}-\deg_\L\Psi_{[1]})
\\
&=\deg\pp_-\Psi_{[1]}-\deg\pp_+\Psi_{[1]}\ge2.
\end{align*}

Let us estimate now the amount of disturbing numbers for the unit disc. Since
the vector functions $\xi$ and $\vt^{\rm t}$ are analytic in $\dd$, the equality
$$
\Psi_{[t_1]}-\Psi_{[t_2]}=(t_1-t_2)u_1\xi\vt^{\rm t}
$$
implies that the multiplicity of the pole at a point $\l\in\dd$ can drop for some $t$
only if $\l$ is a pole of $u_1$. Definitely, each pole of $u_1$ can produce at most one disturbing number.
Thus the total amount of disturbing numbers does not exceed the number of poles of $u_1$
in $\dd$, what is no more than $\deg\pp_-u_1$. \hfill$\bl$

\medskip

The following result gives us a sharp estimate for $\deg\pp_+\Psi$ in the general case. The
example given in the next section shows that this estimate is the best possible.

\begin{thm}
\label{vba}
Let $\Psi$ be a very badly approximable rational $2\times2$
matrix function such that $\deg\pp_-\Psi=k>0$. If $t_1(\Psi)=0$, then
$$
\deg\pp_+\Psi\le k-1.
$$
If $t_1(\Psi)\ne0$, then
$$
\deg\pp_+\Psi\le 2k-3.
$$
\end{thm}

\Pf Consider first the case when $t_1(\Psi)=0$. Then $\Psi$ can be represented as
$$
\Psi=\left(\begin{array}{cc}\bs{w}^\#&\xi\end{array}\right)
\left(\begin{array}{cc}u_0&\0\\\0&\0\end{array}\right)
\left(\begin{array}{c}\bs{v}^*\\\vt^{\rm t}\end{array}\right)
=u_0\bs{w}^\#\bs{v}^*,
$$
where $u_0$ is a scalar badly approximable function,
$\bs{v}$ and $\bs{w}$ are inner and co-outer column functions, and $\xi$ and
$\vt$ are defined by \rf{dop}. By \rf{ner},
$$
\deg\pp_+u_0<\deg\pp_-u_0\le k.
$$
Since $\rank u_0\bs{w}^\#\bs{v}^*=1$ and both $\bs{w}^\#$ and
$\bs{v}^*$ have poles only in $\dd$, it follows easily from Lemma \ref{polu} that
$$
\deg\pp_+\Psi=\deg\pp_+u_0\bs{w}^\#\bs{v}^*\le\deg\pp_+u_0\le k-1.
$$

Suppose now that $t_1(\Psi)\ne0$. Consider a nondisturbing number
$t\in(0,1)$. By Theorem~\ref{gen}
$$
\deg\pp_+\Psi_{[t]}\le\deg\pp_-\Psi_{[t]} -2.
$$
By \rf{Phi},
\bay
\label{u1}
\Psi=\Psi_{[t_1(\Psi)]}=\Psi_{[t]}+(t_1(\Psi)-t)u_1\xi\vt^{\rm t}.
\ey
Since $t$ is
nondisturbing, we have
$$
\deg\pp_+\Psi\le\deg\pp_+\Psi_{[t]}\le\deg\pp_-\Psi_{[t]}-2.
$$
On the other hand, it follows from \rf{u1} that
$$
\deg\pp_-\Psi\ge\deg\pp_-\Psi_{[t]}-\deg\pp_-u_1\xi\vt^{\rm
t}\ge\deg\pp_-\Psi_{[t]}-\deg\pp_-u_1
$$
because both $\xi$ and $\vt^{\rm t}$ are analytic in $\dd$. Thus
$$
\deg\pp_+\Psi-\deg\pp_-\Psi\le\deg\pp_-u_1-2.
$$
The desired inequality follows from \rf{D}. \hfill$\bl$

\begin{thm}
\label{soa}
Let $\Phi$ be a rational $2\times2$ matrix function such that $\deg\Phi=k>0$ and
$\pp_-\Phi\ne\0$. If $t_1(\Phi)=0$\textup, then
$$
\deg\A\Phi\le k-1.
$$
If $t_1(\Phi)\ne0$, then
$$
\deg\A\Phi\le 2k-3.
$$
\end{thm}

\Pf It suffices to apply Theorem~\ref{vba} for the matrix function
$\Psi=\Phi-\A\Phi$. Indeed, Theorem~\ref{vba} yields
$$
\deg(\pp_+\Phi-\A\Phi)\le
\begin{cases}
\phantom2\deg\pp_-\Phi-1\qquad\text{if}\quad t_1(\Phi)=0,
\\
2\deg\pp_-\Phi-3\qquad\text{if}\quad t_1(\Phi)\ne0.
\end{cases}
$$
Substituting this in the inequality $\deg\A\Phi\le\deg\pp_+\Phi+\deg(\pp_+\Phi-\A\Phi)$, we get
$$
\deg\A\Phi\le
\begin{cases}
\deg\pp_+\Phi+\deg\pp_-\Phi-1=\deg\Phi-1\qquad&\text{if}\quad t_1(\Phi)=0,
\\
\deg\pp_+\Phi+2\deg\pp_-\Phi-3\le2\deg\Phi-3\qquad&\text{if}\quad t_1(\Phi)\ne0.
\end{cases}
$$
\hfill$\bl$

\medskip

Theorem \ref{vba} shows that for every $t\in(0,1)$,
$$
\deg\pp_+\Psi_{[t]}+2-\deg\pp_-\Psi_{[t]}\le\deg\pp_-u_1.
$$

The next theorem improves this result. It shows that on the left-hand side of the above inequality we can take 
the sum over all $t\in(0,1]$.

\begin{thm}
\label{an}
Let $\{\Psi_{[t]}\}$, $t\in(0,1]$, be the family of very badly approximable matrix functions defined
by {\em\rf{Phit}}. Then
\bay
\label{mp}
\sum_{t\in(0,1]}\big[\deg\pp_+\Psi_{[t]}+2-\deg\pp_-\Psi_{[t]}\big]_+\le\deg\pp_-u_1.
\ey
\end{thm}

Recall that $[x]_+=\max\{x,0\}$.

\medskip

\Pf By Theorem \ref{gen}, only finite set of disturbing points $t_1,t_2,\cdots,t_n\in(0,1)$,
$n\le\deg\pp_-u_1$, for the unit disk can produce a nonzero term in \rf{mp}.
Denote by
$\L_j$ the subset of $\dd$ such that $t_j$ is a nondisturbing number for $\dd\setminus\L_j$
and it is a disturbing number for every point $\l$ in $\L_j$. Take any nondisturbing
number
$t$ and use the same argument as in the proof of Theorem~\ref{vba}. The only difference
is the following: we estimate the difference of the degrees of $\Psi_{[t_j]}$ and $\Psi_{[t]}$
not for the whole unit disk, but for each $\L_j$ separately. The identity
$$
\Psi_{[t_j]}=\Psi_{[t]}+(t_j-t)u_1\xi\vt^{\rm t}
$$
implies that
$$
\deg_{\L_j}\Psi_{[t_j]}\ge\deg_{\L_j}\Psi_{[t]}-\deg_{\L_j}u_1.
$$
Since $t_j$ is nondisturbing for $\dd\setminus\L_j$, we have the following identity for this set:
$$
\deg_{\dd\setminus\L_j}\Psi_{[t_j]}=\deg_{\dd\setminus\L_j}\Psi_{[t]}.
$$
Thus together with the previous inequality, we obtain
$$
\deg\pp_-\Psi_{[t_j]}\ge\deg\pp_-\Psi_{[t]}-\deg_{\L_j}u_1.
$$
Since $t$ is nondisturbing, we have as before:
$$
\deg\pp_+\Psi_{[t_j]}\le\deg\pp_+\Psi_{[t]}\le\deg\pp_-\Psi_{[t]}-2.
$$
Therefore
$$
\deg\pp_+\Psi_{[t_j]}+2-\deg\pp_-\Psi_{[t_j]}\le\deg_{\L_j}u_1.
$$
Taking the sum over $j$, we obtain \rf{mp}. \hfill$\bl$
\medskip

We will see in Theorem \ref{ekp} that inequality \rf{mp} is best possible
in the sense that if $\vk$ is an integer-valued nonnegative function on
$(0,1]$ that can be nonzero only at finitely many points, then there is a family
of rational very badly approximable functions $\Psi_{[t]}$ of the form \rf{Phit}
such that 
$$
\big[\deg\pp_+\Psi_{[t]}+2-\deg\pp_-\Psi_{[t]}\big]_+=\vk(t)
$$
and
$$
\sum_{t\in(0,1]}\vk(t)=\deg\pp_-u_1.
$$

\

\section{\bf Counterexamples}
\setcounter{equation}{0}

\

In this section for $k\ge2$, we construct an example of a very badly approximable
$2\times2$ rational matrix function
$\Psi$ such that $\deg\pp_-\Psi=k$, but $\deg\pp_+\Psi=2k-3$.
We are going to look for $\Psi$ in the form \rf{Phit}. It follows from
Theorem \ref{gen} that $t$ must be a disturbing number for the unit disk.

We also construct an example of the family $\Psi_{[t]}$ of the form \rf{Phit} such that
inequality \rf{nd} can be violated for several values of $t\in(0,1]$. Actually,
it can be violated for precisely $\deg\pp_-u_1$ values of $t$ which shows that Theorem
\ref{gen} cannot be improved. Moreover, the same construction shows that the inequality in 
Theorem \ref{an} cannot be improved.

The following result
gives a necessary condition for the existence of a disturbing number for the unit disk that
reduces  the multiplicity of a pole at a given point of $\dd$. Its proof is similar to the proof
of Lemma \ref{krat}.

\begin{lem}
\label{vkr}
Suppose that a point $\l\in\dd$ is a pole of $u_1$ of multiplicity $d>0$
and the multiplicity of the pole of $\Psi$  at $\l$ is $d-l$, where $l>0$.
Then $u_0$ has a zero at $\l$ of multiplicity at least $l$ and both $\bs{v}^\#$
and $\bs{w}^\#$ have a
pole at $\l$ of multiplicity at least $l$.
\end{lem}

\Pf It is easy to see that
\bay
\label{w}
\bs{w}^{\rm t}\Psi=u_0\bs{v}^*
\ey
and
\bay
\label{v}
\Psi\bs{v}=u_0\bs{w}^\#.
\ey

Since by Lemma \ref{dt}, $\deg_\l\det\Psi\le d-l$, it follows from \rf{det} that $u_0$
has a zero at $\l$ of multiplicity at least $l$.

Multiplying~\eqref{w} on the right and~\eqref{v} on the left
by the adjugate matrix $\Psi^{\rm ad}$ and using~\eqref{det}, we get

\bay 
\label{w1} 
tu_1\bs{w}^{\rm t}=\bs{v}^*\Psi^{\rm ad} \ey and \bay \label{v1}
tu_1\bs{v}=\Psi^{\rm ad}\bs{w}^\#. 
\ey 
By Lemma \ref{jug},
$\deg_\l\Psi^{\rm ad}=\deg_\l\Psi$. Since $\bs{v}(\l)\ne0$ and $\bs{w}^{\rm t}(\l)\ne0$
(their components are coprime!), the expressions in these equalities have a pole of
multiplicity $d$ at the point $\l$. Therefore both $\bs{v}^*$ (and hence $\bs{v}^\#$) and
$\bs{w}^\#$ have poles at $\l$, of multiplicity at least $l$. \hfill$\bl$

\medskip

We can try to construct a desired example of a function $\Psi$, for which all poles
have multiplicity one. To construct such an example, we have to cancel $k-1$ poles of
$tu_1\xi\vt^{\rm t}$ in $\dd$. Note that $\deg\pp_-u_1\le k-1$ (see \S\,2.3). Thus we have to
cancel all poles of $tu_1\xi\vt^{\rm t}$ in $\dd$. By Lemma \ref{vkr}, $u_0$ must have zeros
at all poles of $u_1$ in $\dd$ and both $\bs{v}^\#$ and $\bs{w}^\#$ must have poles at all
poles of $u_1$ in $\dd$.

\begin{thm}
\label{kp} Let $k\ge2$ be an integer\textup, $0<t<1$, and $a>1$.
Suppose that $B_0$\textup, $B_1$\textup, and $B_2$ are Blaschke products  such that
$$
\deg B_0=k-1,\quad\deg B_1=k,\quad\deg B_2=k-2,
$$
and the function $1-ta^2B_1B_2$ vanishes at the zeros of $B_0$. Let
\bay
\label{psi} \Psi=
\left(\begin{array}{cc}\frac1{B_0}&a\\-a&B_0\end{array}\right)
\left(\begin{array}{cc}\frac{B_0}{B_1}&\0\\\0&t\frac{B_2}{B_0}\end{array}\right)
\left(\begin{array}{cc}\frac1{B_0}&a\\-a&B_0\end{array}\right).
\ey
Then
$\deg\pp_-\Psi=k$\textup, while $\deg\pp_+\Psi=2k-3$.
\end{thm}

{\bf Remark.} On the right-hand side of formula \rf{psi} the first and the third
matrix functions are not as in the definition of thematic factorizations (see \rf{si}).
However, if we multiply the first and the third matrix functions on the right-hand side of
\rf{psi} by $(1+a^2)^{1/2}$, we obtain matrix functions as in the thematic factorization
\rf{si}. Thus the function $\Psi$ defined by \rf{psi} is very badly approximable.

\medskip

\Pf It is easy to verify that
$$
\Psi=\left(\begin{array}{cc}\frac{1-ta^2B_1B_2}{B_0B_1}&\frac{a}{B_1}+taB_2\\[.2cm]
-\frac{a}{B_1}-taB_2&B_0(\frac{-a^2}{B_1}+tB_2)
\end{array}\right).
$$

Let us first compute $\deg\pp_-\Psi$. By the hypotheses of the theorem, $\pp_-\Psi$ has removable
singularities at the zeros of $B_0$. Thus $\pp_-\Psi$ has poles only at the zeros of $B_1$.
Hence,
$$
\deg\pp_-\Psi=
\deg_{B_1^{-1}(0)}\left(\begin{array}{cc}\frac1{B_0B_1}&\frac{a}{B_1}\\[.2cm]
-\frac{a}{B_1}&-a^2\frac{B_0}{B_1}
\end{array}\right).
$$
We have
$$
\left(\begin{array}{cc}\frac1{B_0B_1}&\frac{a}{B_1}\\[.2cm]
-\frac{a}{B_1}&-a^2\frac{B_0}{B_1}
\end{array}\right)=\frac{1}{B_1}\left(\begin{array}{cc}\frac1{B_0}&a\\[.2cm]
-a&-a^2 B_0
\end{array}\right).
$$
Clearly,
$$
\rank\left(\begin{array}{cc}\frac1{B_0}&a\\[.2cm]
-a&-a^2 B_0
\end{array}\right)=1,
$$
and so
$$
\deg\pp_-\Psi=\deg B_1=k.
$$

Let us now compute $\deg\pp_+\Psi$. Clearly, $\pp_+\Psi$ can have poles only at the poles of
$B_2$ and the poles of $B_0$.

Let us first show that there can be no cancellation of poles at the poles of
$B_2$ or $B_0$. The only possibility of cancellation is to cancel poles of $B_0$
of the lower right entry which can be done if $-\frac{a^2}{B_1}+tB_2$ vanishes at poles of $B_0$.
Let $\z$ be a zero of $B_0$ and suppose that
$$
\frac{a^2}{B_1(1/\bar\z)}=tB_2(1/\bar\z).
$$
It follows that
$$
a^2B_1(\z)B_2(\z)=t.
$$
Since $t<1$, this contradicts the assumption that $ta^2B_1(\z)B_2(\z)=1$.

First, we compute the contribution of the poles of $B_0$. Clearly, only the lower right entry
of $\pp_+\Psi$ has poles at poles of $B_0$. Obviously,
$$
\deg_{B_0^{-1}(\be)}\Psi=k-1.
$$
Let us compute the contribution of the poles of $B_2$. We have
$$
\deg_{B_2^{-1}(\be)}\Psi=
\deg_{B_2^{-1}(\be)}\left(\begin{array}{cc}\frac{-ta^2B_2}{B_0}&taB_2\\[.2cm]
-taB_2&tB_2B_0
\end{array}\right)
=\deg_{B_2^{-1}(\be)}B_2\left(\begin{array}{cc}\frac{-ta^2}{B_0}&ta\\[.2cm]
-ta&tB_0\end{array}\right).
$$
Clearly,
$$
\rank\left(\begin{array}{cc}\frac{-ta^2}{B_0}&ta\\[.2cm]
-ta&tB_0
\end{array}\right)=1,
$$
and so the contribution of the poles of $B_2$ is $\deg B_2=k-2$. Thus
$$
\quad\qquad\qquad\qquad\qquad\qquad\deg\pp_+\Psi=k-1+k-2=2k-3.\quad\qquad\qquad\qquad\qquad\qquad\bl
$$

\begin{cor}
\label{sp} Let $\Psi$ be the function defined in Theorem~\textup{\ref{kp}} and let
$\Phi=\pp_-\Psi$. Then $\deg\Phi=k$ and $\deg\A\Phi=2k-3$.
\end{cor}

Finally, we can slightly modify the construction in Theorem \ref{kp} to get a family
$\Psi_{[t]}$ of the form \rf{Phit} such that inequality \rf{nd} is violated
for several different values of $t\in(0,1)$.

The following result shows that roughly speaking, inequality \rf{mp} is the only constraint
for the numbers $\big[\deg\pp_+\Psi_{[t]}+2-\deg\pp_-\Psi_{[t]}\big]_+$.

\begin{thm}
\label{ekp}
Let $k$\textup, $a$\textup, $B_0$\textup, $B_1$\textup, and $B_2$ be as in
Theorem~\textup{\ref{kp}}. Suppose that $t_1,\cdots,t_d$ are distinct points in
$(0,1)$\textup, $\D_1,\cdots,\D_d$ is a partition of the zero set of $B_0$\textup, and
$\vk_j$ is the number of zeros of $B_0$ in $\D_j$ (counted with multiplicities). Assume that
\bay
\label{int}
1-t_ja^2B_1B_2\quad
\mbox{vanishes at the points of}\quad\D_j,\quad1\le j\le d.
\ey
If the very badly
approximable function $\Psi_{[t]}$ is defined by
$$
\Psi_{[t]}=
\left(\begin{array}{cc}\frac1{B_0}&a\\-a&B_0\end{array}\right)
\left(\begin{array}{cc}\frac{B_0}{B_1}&\0\\\0&t\frac{B_2}{B_0}\end{array}\right)
\left(\begin{array}{cc}\frac1{B_0}&a\\-a&B_0\end{array}\right),\quad 0<t\le1,
$$
then inequality {\em\rf{nd}} is violated precisely for $t=t_1,\cdots,t_d$
and
$$
\deg\pp_+\Psi_{[t_j]}=\deg\pp_-\Psi_{[t_j]}-2+\vk_j.
$$
\end{thm}

The proof of Theorem \ref{ekp} is exactly the same as the proof of Theorem \ref{kp}.

In particular, if each $\D_j$ is a singleton, inequality \rf{nd} is violated for precisely
$\deg\pp_-u_1$ values of $t$.

\medskip

{\bf Remark.} Let us fix the zeros of $B_0$. Condition \rf{int} means that the Blaschke
product $B_1B_2$ interpolates at the points of $\D_j$ the value $(a^2t_j)^{-1}$. It is well
known (see e.g., \cite{Pe}, Ch.~7, \S\,1) that if $a$ is sufficiently large, it is possible to
find a Blaschke product of a given degree greater than or equal to $k$ that solves this
interpolation problem. Thus if $a$ is sufficiently large, it is possible to find $B_1$ and
$B_2$ satisfying the hypotheses of the theorem.

\

\

\noindent
\begin{tabular}{p{8cm}p{14cm}}
V.V. Peller & V.I. Vasyunin \\
Department of Mathematics & St.Petersburg Branch \\
Michigan State University  & Steklov Institute of Mathematics \\
East Lansing, Michigan 48824 & Fontanka 27, 191011 St-Petersburg\\
USA&Russia
\end{tabular}

\end{document}